\documentclass[12pt]{article}

\usepackage{amsfonts,amsmath,amssymb}
\usepackage{mathrsfs,mathtools,url}
\usepackage{latexsym}
\usepackage{tikz,color}

\newtheorem{theorem}{Theorem}

\newtheorem{lemma}[theorem]{Lemma}
\newtheorem{proposition}[theorem]{Proposition}

\newtheorem{remark}[theorem]{Remark}

\DeclareMathOperator {\diam} {diam}

\DeclareMathOperator {\mut} {\mu_{\rm t}}

\def\cp{\,\square\,}
\newcommand{\proof}{\noindent{\bf Proof.\ }}
\newcommand{\qed}{\hfill $\square$ \bigskip}
\newcommand{\smallqed}{{\tiny ($\Box$)}}
\newcommand{\ex}{{\rm ex}}
\newcommand{\cK}{{\cal K}}
\newcommand{\cO}{{\cal O}}

\vspace{8mm}
\title{Total mutual-visibility in Hamming graphs}

\author{Csilla Bujt\'{a}s$^{a, b}$, Sandi Klav\v{z}ar$^{a, b, c}$, Jing Tian$\/^{d, e}$ \\\\
$^{a}$ \small Faculty of Mathematics and Physics, University of Ljubljana, Slovenia\\
\small {\tt csilla.bujtas@fmf.uni-lj.si}\\
\small {\tt sandi.klavzar@fmf.uni-lj.si}\\
$^{b}$ \small Institute of Mathematics, Physics and Mechanics, Ljubljana, Slovenia \\
$^{c}$ \small Faculty of Natural Sciences and Mathematics, University of Maribor, Slovenia\\
$^{d}$ \small School of Science, Zhejiang University of Science and Technology, \\
\small Hangzhou, Zhejiang 310023, PR China\\
\small {\tt jingtian526@126.com}\\
}
\date{}

\begin{document}

\maketitle

\begin{abstract}
If $G$ is a graph and $X\subseteq V(G)$, then $X$ is a total mutual-visibility set if every pair of vertices $x$ and $y$ of $G$ admits a shortest $x,y$-path $P$ with $V(P) \cap X \subseteq \{x,y\}$. The cardinality of a largest total mutual-visibility set of $G$ is the total mutual-visibility number $\mu_{\rm t}(G)$ of $G$. In this paper the total mutual-visibility number is studied on Hamming graphs, that is, Cartesian products of complete graphs. Different equivalent formulations for the problem are derived. The values $\mu_{\rm t}(K_{n_1}\,\square\, K_{n_2}\,\square\, K_{n_3})$ are determined. It is proved that $\mu_{\rm t}(K_{n_1} \,\square\, \cdots \,\square\, K_{n_r}) = \cO(N^{r-2})$, where $N = n_1+\cdots + n_r$, and that $\mu_{\rm t}(K_s^{\,\square\,, r}) = \Theta(s^{r-2})$ for every $r\ge 3$, where $K_s^{\,\square\,, r}$ denotes the Cartesian product of $r$ copies of $K_s$. The main theorems are also reformulated as Tur\'an-type results on hypergraphs.
\end{abstract}

\noindent
\textbf{Keywords:} mutual-visibility set; total mutual-visibility set; Hamming graph; Tur\'an-type problem

\medskip\noindent
\textbf{AMS Math.\ Subj.\ Class.\ (2020)}: 05C12, 05C38, 05C65, 05C76


\section{Introduction}

Let $G = (V(G), E(G))$ be a graph and $X\subseteq V(G)$. Then vertices $x$ and $y$ of $G$ are {\em $X$-visible}, if there exists a shortest $x,y$-path $P$ such that no internal vertex of $P$ belongs to $X$. The set $X$ is a \emph{mutual-visibility set} if any two vertices from $X$ are $X$-visible, while $X$ is a \emph{total mutual-visibility set} if any two vertices from $V(G)$ are $X$-visible. The cardinality of a largest  mutual-visibility set (resp.\ total mutual-visibility set) is the \emph{mutual-visibility number} (resp.\ \emph{total mutual-visibility number}) $\mu(G)$ (resp.\ $\mut(G)$) of $G$.

The mutual-visibility sets were introduced by Di Stefano in~\cite{distefano-2022} motivated by mutual visibility in distributed computing and social networks. Although the motivation came from theoretical computer science, it is a graph theory concept. It needs to be said that the term mutual-visibility is also used in other contexts, for instance in robotics, where the mutual visibility problem asks for a distributed algorithm that reposition robots to a configuration where they all can see each other, cf.~\cite{adhikary-2022}. Some related research can be found in~\cite{Cicerone-2023+, diluna-2017, poudel-2021}. The graph theoretic mutual-visibility was further investigated in~\cite{cicerone-2023, cicerone-2023+}, where the latter paper naturally raised the need to introduce the total mutual-visibility which was in turn investigated in~\cite{kuziak-2023+, tian-2023+}.

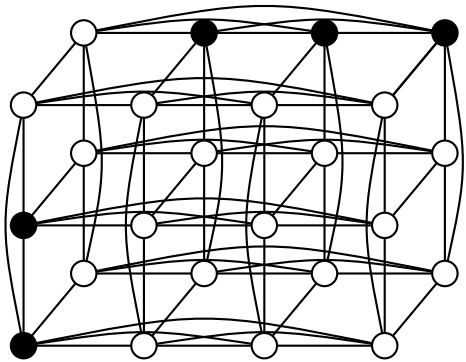
\begin{figure}[ht!]
\begin{center}
\begin{tikzpicture}[scale=1.6,style=thick]
\tikzstyle{every node}=[draw=none,fill=none]
\def\vr{3pt} 

\begin{scope}[yshift = 0cm, xshift = 0cm]
\path (0,0) coordinate (x1);
\path (1,0) coordinate (x2);
\path (2,0) coordinate (x3);
\path (3,0) coordinate (x4);
\path (3.5,0.6) coordinate (x5);
\path (2.5,0.6) coordinate (x6);
\path (1.5,0.6) coordinate (x7);
\path (0.5,0.6) coordinate (x8);
\path (0,1) coordinate (x9);
\path (1,1) coordinate (x10);
\path (2,1) coordinate (x11);
\path (3,1) coordinate (x12);
\path (3.5,1.6) coordinate (x13);
\path (2.5,1.6) coordinate (x14);
\path (1.5,1.6) coordinate (x15);
\path (0.5,1.6) coordinate (x16);
\path (0,2) coordinate (x17);
\path (1,2) coordinate (x18);
\path (2,2) coordinate (x19);
\path (3,2) coordinate (x20);
\path (3.5,2.6) coordinate (x21);
\path (2.5,2.6) coordinate (x22);
\path (1.5,2.6) coordinate (x23);
\path (0.5,2.6) coordinate (x24);

\draw (x1) --(x2)--(x3)--(x4)--(x5)--(x6)--(x7)--(x8);
\draw (x9)--(x10)--(x11)--(x12)--(x13)--(x14)--(x15)--(x16)--(x9);
\draw (x17)--(x18)--(x19)--(x20)--(x21)--(x22)--(x23)--(x24);
\draw (x1) -- (x9)--(x17) -- (x24)--(x16) -- (x8)--(x1);
\draw (x2) -- (x10)--(x18) -- (x23)--(x15) -- (x7)--(x2);
\draw (x3) -- (x11)--(x19) -- (x22)--(x14) -- (x6)--(x3);
\draw (x4) -- (x12)--(x20) -- (x21)--(x13) -- (x5);
\draw (x10) -- (x15);
\draw (x11) -- (x14);
\draw (x1) .. controls (1,0.15).. (x3);
\draw (x2) .. controls (2,0.15).. (x4);
\draw (x1) .. controls (1.5,0.3).. (x4);

\draw (x8) .. controls (1.5,0.75).. (x6);
\draw (x7) .. controls (2.5,0.75).. (x5);
\draw (x8) .. controls (2,0.9).. (x5);

\draw (x9) .. controls (1,1.15).. (x11);
\draw (x10) .. controls (2,1.15).. (x12);
\draw (x9) .. controls (1.5,1.3).. (x12);
\draw (x16) .. controls (1.5,1.75).. (x14);

\draw (x15) .. controls (2.5,1.75).. (x13);
\draw (x16) .. controls (2,1.9).. (x13);

\draw (x17) .. controls (1,2.15).. (x19);
\draw (x18) .. controls (2,2.15).. (x20);
\draw (x17) .. controls (1.5,2.3).. (x20);
\draw (x24) .. controls (1.5,2.75).. (x22);
\draw (x23) .. controls (2.5,2.75).. (x21);
\draw (x24) .. controls (2,2.9).. (x21);

\draw (x1) .. controls (-0.2,1).. (x17);
\draw (x8) .. controls (0.7,1.5).. (x24);
\draw (x2) .. controls (0.8,1).. (x18);
\draw (x7) .. controls (1.7,1.5).. (x23);
\draw (x3) .. controls (1.8,1).. (x19);
\draw (x6) .. controls (2.7,1.5).. (x22);
\draw (x4) .. controls (2.8,1).. (x20);
\draw (x5) .. controls (3.7,1.5).. (x21);


\draw (x1)  [fill=black] circle (\vr);
\draw (x2)  [fill=white] circle (\vr);
\draw (x3)  [fill=white] circle (\vr);
\draw (x4)  [fill=white] circle (\vr);
\draw (x5)  [fill=white] circle (\vr);
\draw (x6)  [fill=white] circle (\vr);
\draw (x7)  [fill=white] circle (\vr);
\draw (x8)  [fill=white] circle (\vr);
\draw (x9)  [fill=black] circle (\vr);
\draw (x10) [fill=white] circle (\vr);
\draw (x11) [fill=white] circle (\vr);
\draw (x12) [fill=white] circle (\vr);
\draw (x13) [fill=white] circle (\vr);
\draw (x14) [fill=white] circle (\vr);
\draw (x15) [fill=white] circle (\vr);
\draw (x16) [fill=white] circle (\vr);
\draw (x17) [fill=white] circle (\vr);
\draw (x18) [fill=white] circle (\vr);
\draw (x19) [fill=white] circle (\vr);
\draw (x20) [fill=white] circle (\vr);
\draw (x21) [fill=black] circle (\vr);
\draw (x22) [fill=black] circle (\vr);
\draw (x23) [fill=black] circle (\vr);
\draw (x24) [fill=white] circle (\vr);

\end{scope}
\end{tikzpicture}
\end{center}
\caption{$K_2\cp K_3\cp K_4$ with a total mutual-visibility set of cardinality $5$ in bold.}
\label{fig:graph 1}
\end{figure}

A graph $G$ is a {\em Hamming graph} if $G$ is the Cartesian product of complete graphs. In particular, complete graphs are Hamming graphs. In~\cite[Corollary 3.7]{cicerone-2023} it was shown that $\mu(K_n\cp K_m) = z(n,m;2,2)$, where $z(n,m;2,2)$ is the Zarankiewitz's number. To determine the latter number is a notorious open problem~\cite{west-2021, z-1951}. On the other hand, it was proved in~\cite[Proposition 15]{tian-2023+} that $\mut(K_n\cp K_m) = \max \{n,m\}$. In~\cite{kuziak-2023+} the authors proposed a challenging problem to determine the total mutual-visibility number of Hamming graphs with at least three factors. They provided a total mutual-visibility set of $K_3\cp K_3\cp K_2$ of cardinality $4$, and in Fig.~\ref{fig:graph 1} we give a total mutual-visibility set of $K_2\cp K_3\cp K_4$ of cardinality $5$.

In the light of what has just been said, in this paper we focus on the total mutual-visibility in Hamming graphs. In the next section we give equivalent formulations of the problem which later serve as a tool for proof. Our first main result which we prove in Section~\ref{sec:3-dim} reads as follows.

\begin{theorem}
\label{thm:three-complete graphs}
If $n_1\geq n_2\geq n_3\ge 2$ and $N = n_1 + n_2 + n_3$, then
\[
  \mut(K_{n_1}\cp K_{n_2}\cp K_{n_3}) = \left.
  \begin{cases}
    N - 4; & n_3 = 2, \\
    N - 5; & n_3 = 3, \\
    N - 6; & n_3\ge 4\,.
  \end{cases}
  \right.
\]
\end{theorem}

Note that if $n_1\geq n_2\geq n_3 = 1$, then by the above-mentioned result from~\cite{tian-2023+} we have $  \mut(K_{n_1}\cp K_{n_2}\cp K_1) = \mut(K_{n_1}\cp K_{n_2}) = N - n_2 -1 = n_1$.

In Section~\ref{sec:upper}, we prove the following:

\begin{theorem}  \label{thm:upper}
If $r \ge 3$, and $s, n_1,\ldots, n_r$ are positive integers, $N = n_1+\cdots + n_r$, then the following statements hold:
\begin{itemize}
\item[$(i)$] $\mut(K_{n_1} \cp \cdots \cp K_{n_r}) \le \frac{6}{r!}\,N^{r-2}$;
\item[$(ii)$] $\mut(K_s^{\cp, r}) \leq c_r'\,s^{r-2}$ with $c_r' =  3\prod_{i=3}^r  (i-1)^{i-3}$. 
\end{itemize}
\end{theorem}

In the subsequent section, we strengthen the result for balanced Hamming graphs by establishing the exact magnitude for their total mutual-visibility number. The result~\cite[Proposition 15]{tian-2023+} implies $\mut( K_s^{\cp, 2} )= \Theta(s)$.  
However, the situation is different for higher values of $r$ as our third main result asserts. It will be proved in Section~\ref{sec:balanced} using a probabilistic approach.

\begin{theorem} \label{thm:lower}
If $r\ge 3$, then
$$\mut(K_s^{\cp, r}) = \Theta(s^{r-2}).$$
\end{theorem}

In the last section we reformulate our problem as a Tur\'an-type question and accordingly restate Theorems~\ref{thm:three-complete graphs}-\ref{thm:lower}.

In the remainder of the introduction, we recall some definitions and terminology, mainly about the Cartesian product of graphs. The standard shortest-path distance between vertices $u$ and $v$ of a (connected) graph $G$ will be denoted by $d_G(u,v)$. We will use the term clique for a complete graph as well as for its vertex set. If $u,v\in V(G)$, then the {\em interval} $I_G[u,v]$ between $u$ and $v$ in $G$ is the set of all vertices of $G$ that lie on shortest $u,v$-paths. The {\em Cartesian product} $G\cp H$ of graphs $G$ and $H$ has the vertex
set $V(G\cp H) = V(G)\times V(H)$, vertices $(g, h$) and $(g', h')$ are adjacent if either $gg'\in E(G)$ and $h = h'$, or $g = g'$ and $hh'\in E(H)$. Given a vertex $h\in V(H)$, the subgraph of $G\cp H$ induced by the set $\{(g,h): g\in V(G)\}$, is a {\em $G$-layer} and is denoted by $G^h$. $H$-layers $^gH$ are defined analogously. Each $G$-layer and each $H$-layer is isomorphic to $G$ and $H$, respectively. Moreover, each layer of a Cartesian product is its convex subgraph. We use this fact later on many times, sometimes implicitly. The Cartesian product of $r$ copies of $G$ is denoted by $G^{\cp, r}$. We say that $K_s^{\cp, r}$ is a \emph{balanced Hamming graph}. For more information on the Cartesian product see the book~\cite{hammack-2011}.

\section{Equivalent formulations of the problem}
\label{sec:equivalent}

In this section we prove two equivalent formulations of the total mutual-visibility problem in Hamming graphs. First we prove that total mutual-visibility sets in Hamming graphs are precisely the vertex sets such that no pair of vertices is at distance $2$. Then we reformulate this fact in terms of clique systems in complete multipartite graphs.

We say that a $4$-cycle of a Cartesian product graph $G$ is a {\em Cartesian square} if it is not contained in a single layer of $G$. This definition also applies to Cartesian product of more than two factors. More precisely, let $G = G_1\cp\cdots \cp G_r$, $r\ge 2$. Then the vertices $u, u', u'', u'''\in V(G)$ form a Cartesian square if there exist $i,j\in [r]$, $i < j$, such that
\begin{align*}
u & = (u_1,\ldots, u_{i-1},u_i,u_{i+1},\ldots, u_{j-1},u_j,u_{j+1},\ldots, u_r)\,, \\
u' & = (u_1,\ldots, u_{i-1},u_i,u_{i+1},\ldots, u_{j-1},u_j',u_{j+1},\ldots, u_r)\,, \\
u'' & = (u_1,\ldots, u_{i-1},u_i',u_{i+1},\ldots, u_{j-1},u_j',u_{j+1},\ldots, u_r)\,, \\
u''' & = (u_1,\ldots, u_{i-1},u_i',u_{i+1},\ldots, u_{j-1},u_j,u_{j+1},\ldots, u_r)\,.
\end{align*}
The following results were stated in~\cite[Lemma~5.8]{kuziak-2023+} for two factors, the proof for more factors is analogous. That is, we just need to infer that a Cartesian square is a convex subgraph of a Cartesian product graph.

\begin{lemma}
\label{lem:basic-cp}
Let $G = H_1 \cp \cdots \cp H_r$, $r\ge 2$. If $X$ is a total mutual-visibility set of $G$ and $C$ is a Cartesian square of $G$, then $X$ contains no diametral pair of vertices of $C$.
\end{lemma}

If a Cartesian square of a Cartesian product graph $G$ fulfils the condition of Lemma~\ref{lem:basic-cp} for a set $X\subseteq V(G)$, then we say that the cycle is {\em $X$-suitable}. For the proof of our first characterization, we need the following result.

\begin{proposition}
\label{prop:intervals-in-Hamming}
If $G$ is a Hamming graph, $u,v\in V(G)$, and $d_G(u,v) = t$, then the subgraph induced by $I[u,v]$ is isomorphic to the $t$-cube $Q_t$.
\end{proposition}

Our first equivalent formulation of the total mutual-visibility problem in Hamming graphs now reads as follows.

\begin{theorem}
\label{thm:Hamming}
If $G$ is a Hamming graph and $X\subseteq V(G)$, then the following statements are equivalent.
\begin{enumerate}
\item[(i)] $X$  is a total mutual-visibility set of $G$.
\item[(ii)] $2\notin \{ d_G(u,v):\ u,v\in X\}$.
\item[(iii)] Each Cartesian square of $G$ is $X$-suitable.
\end{enumerate}
\end{theorem}

\proof
Let $G = K_{n_1}\cp \cdots \cp K_{n_r}$, where $n_i\ge 2$ for $i\in [r]$, and $r\ge 1$. If $r=1$, then $G$ is a complete graph which contains no Cartesian square, every subset of $V(G)$ forms a total mutual-visibility set, and $\diam(G) = 1$. Hence the three statements are equivalent for $G$ and we may assume in the rest that $r\ge 2$.

$(i) \Rightarrow (ii)$: Suppose on the contrary that there exist vertices $u,v\in X$ with $d_G(u,v) = 2$. Then $u$ and $v$ lie in a convex $C_4$, but then the other two vertices of this convex $C_4$ are not $X$-visible, a contradiction.

$(ii) \Rightarrow (iii)$: If $2\notin \{ d_G(u,v):\ u,v\in X\}$, then clearly each Cartesian square contains at most two vertices of $X$, and if so, these two vertices are adjacent, hence each Cartesian square is $X$-suitable.

$(iii) \Rightarrow (i)$:
Let $X\subseteq V(G)$ and assume that each Cartesian square of $G$ is $X$-suitable. We need to show that each two vertices $u,v\in V(G)$ are $X$-visible and proceed by induction on $d_G(u,v) = t$. If $t=1$, the assertion is clear. Suppose now that $t\ge 2$ and that each pair of vertices at distance at most $t-1$ is $X$-visible. By Proposition~\ref{prop:intervals-in-Hamming}, $I_G[u,v]$ induces a $t$-cube $Q_t$. Let $v^{(1)}, \ldots, v^{(t)}$ be the neighbors of $v$ in this $Q_t$. Then each pair of vertices $v^{(i)}$ and $v^{(j)}$ lies in a Cartesian square of $G$ together with the vertex $v$. Since each Cartesian square of $G$ is $X$-suitable, we infer that at most one of the vertices $v^{(1)}$ and $v^{(2)}$ lies in $X$. Assume without loss of generality that $v^{(1)}\notin X$. By the induction hypothesis, $v^{(1)}$ and $u$ are $X$-visible which in turn implies that $u$ and $v$ are also $X$-visible by concatenating the corresponding shortest $u,v^{(1)}$-path and the edge $v^{(1)}v$.
\qed

To get another reformulation of the total mutual-visibility problem on Hamming graphs, we consider the complete multipartite graph $K_{n_1,\ldots,n_r}$, where $r\ge 3$ and $n_1\ge \cdots \ge n_r\ge 2$. Note that each maximal clique in it is a maximum clique that is of order $r$.

\begin{proposition}
	\label{prop:reduction}
If $r\ge 3$ and $n_1\ge \cdots \ge n_r\ge 2$, then $\mut(K_{n_1}\cp  \cdots\cp K_{n_r})$ is equal to the cardinality of a largest family of maximal cliques of $K_{n_1,\ldots,n_r}$ such that no two cliques from the family intersect in an $(r-2)$-clique.
\end{proposition}

\proof
Let $r\ge 3$ and $n_1\ge \cdots \ge n_r\ge 2$ and set $G = K_{n_1}\cp \cdots\cp K_{n_r}$. Setting $V(K_{n_j})=[n_j]$ we have $V(G) = \{(i_1,\ldots, i_r):\ i_j\in [n_j], j\in [r]\}$. We are going to prove that each total mutual-visibility set of $G$ gives rise to a family of maximal cliques in $K_{n_1,\ldots,n_r}$, such that no two cliques from the family intersect in an $(r-2)$-clique, as well as the other way around, that is, each family of maximal cliques in $K_{n_1,\ldots,n_r}$ gives rise to a total mutual-visibility set of $G$.

Set $H = K_{n_1,\ldots,n_r}$. Let $I_i$, $i\in [r]$, be the partite classes of $H$, where $|I_i| = n_i$, so that $V(H) = \bigcup_{i=1}^r I_i$. Let further $I_i = \{u_{i,j}:\ j \in [n_i]\}$.

Let $X = \{x_1,\ldots, x_t\}$ be a total mutual-visibility set of $G$. For $i\in [t]$, set
\begin{equation}
	\label{eq:first}
	x_i = (z_1^{(i)}, z_2^{(i)}, \ldots, z_r^{(i)})\,.
\end{equation}
To each vertex $x_i$ assign an $r$-clique of $H$ induced by the set of vertices
\begin{equation}
	\label{eq:second}
	X_i = \{u_{1,z_1^{(i)}}, u_{2,z_2^{(i)}}, \dots, u_{r,z_r^{(i)}}\}\,.
\end{equation}
We claim that ${\cal X} = \{X_i:\ i\in [t]\}$ is a set of $r$-cliques of $H$ such that no two cliques from ${\cal X}$ intersect in an $(r-2)$-clique. Since the vertices from $X_i$ belong to pairwise different partite classes of $H$, the graph induced by $X_i$ is isomorphic to $K_r$. Moreover, if $i'\ne i$, then $d_G(x_i, x_{i'}) \ne 2$, hence $|X_i \cap X_{i'}| \ne r-2$. We have thus seen that the total mutual-visibility set $X$ of $G$ gives rise to a family of maximal cliques of 
$H$, such that no two cliques from the family intersect in an $(r-2)$-clique.

To prove the reverse assignment, we proceed in the reverse order as above. That is, we start with a family of $k$-cliques ${\cal X}$ such that no two cliques from the set intersect in an $(r-2)$-clique. Then we use their enumeration as in~\eqref{eq:second} to produce a total mutual-visibility set of the same cardinality as in~\eqref{eq:first}.
\qed

\section{Proof of Theorem~\ref{thm:three-complete graphs}}
\label{sec:3-dim}

The proof of Theorem~\ref{thm:three-complete graphs} is divided into two cases. We first deal with the case when $n_3\in \{2,3\}$, and then with the case $n_3\ge 4$. For the proof of the first part, we recall the following result.

\begin{proposition} {\rm \cite[Proposition 4.4]{tian-2023+}}
\label{prop:cp-cycle-by-complete graphs}
If $s\geq 3$ and $n\geq 3$, then
 $$\mut(C_s\cp K_n) = \left\{
\begin{array}{ll}
0; & s\geq 5,\\
n; & \mbox{otherwise}.\\
\end{array}\right.
$$
\end{proposition}

\begin{theorem}
\label{thm:three-complete graphs-small-cases}
If $n_1\geq n_2\geq n_3\in \{2,3\}$, then $\mut(K_{n_1}\cp K_{n_2}\cp K_{n_3}) = n_1+n_2-2$.
\end{theorem}

\proof
Let $n_1\geq n_2\geq n_3$, let $n_3\in \{2,3\}$, and set $G=K_{n_1}\cp K_{n_2}\cp K_{n_3}$. By a straightforward case analysis we infer that $\mut(K_2\cp K_2\cp K_2) = 2$, hence the theorem holds in this case. We may thus assume in the rest that $n_1\ge 3$.

The vertices from the set $Y = \{(i,1,1):\ 2\le i\le n_1\} \cup \{(1,j,2):\ 2\le j\le n_2\}$ are pairwise at distance $1$ or $3$. Theorem~\ref{thm:Hamming} implies that $Y$ is a total mutual-visibility set of $G$. Hence $\mut(G) \ge n_1+n_2-2 = |Y|$.

Let $V(K_{n_j})=[n_j]$, so that $V(G) = \{(i_1,i_2,i_3):\ i_j\in [n_j], j\in [3]\}$ and let $X$ be a total mutual-visibility set with  $|X| = \mut(G)$. For a vertex $x\in V(G)$, let $X_i(x) = \{v\in X:\ d_G(x,v)=i\}$, $i\in [3]$.

To prove that $\mut(G) \le n_1+n_2-2$, note first that if no two vertices of $X$ are adjacent, then $|X| \le n_3$. Indeed, if $|X| > n_3$, then there exist two vertices $w,w'\in X$ with the same third coordinate. As $w$ and $w'$ are not adjacent, this means then $d_G(w,w') = 2$, a contradiction with Theorem~\ref{thm:Hamming}. Hence, if no two vertices of $X$ are adjacent, then $|X| \le n_3 \le n_2 \le n_1 + n_2 - 2$.

Let $u$ and $u'$ be two adjacent vertices of $X$. Assume that $u$ and $u'$ differ in the first coordinate. By the symmetry of Hamming graphs we may assume without loss of generality that $u = (1,1,1)$ and $u' = (2,1,1)$. By Theorem~\ref{thm:Hamming}, we have $X_2(u) = \emptyset$, so that
\begin{equation}
\label{eq:union}
X = \{u\} \cup X_1(u)\cup X_3(u)\,.
\end{equation}
We claim that all the vertices from $X_1(u)$ differ from $u$ in the first coordinate. Indeed, if this would not be the case, then there would exist a vertex $u'' = (1,j'',1)$ (or $u'' = (1,1,k'')$) in $X$, but then $u'$ and $u''$ are diametral vertices of a Cartesian square which is not possible by Theorem~\ref{thm:Hamming}. Since $n_1\ge n_2\ge n_3$, the claim implies that $|X_1|\le n_1-1$.

\medskip\noindent
{\bf Case 1}: $|X_1(u)|= n_1-1$. \\
By the above claim, in this case we may assume without loss of generality that $X_1(u) = \{(i,1,1):\ i\in \{2,\ldots, n_1\}\}$. Consider an arbitrary vertex $w = (w_1,w_2,w_3)$ from $X_3(u)$. Then $w_1\ne 1$,  $w_2\ne 1$, and  $w_3\ne 1$. Since $(w_1,1,1)\in X$ and $d(w, (w_1,1,1)) = 2$ we get that $X_3(u) = \emptyset$.  By~\eqref{eq:union} we conclude that if $|X_1(u)|= n_1-1$, then $|X| = n_1 \le n_1 + n_2 - 2$ because $n_2\ge 2$.

\medskip\noindent
{\bf Case 2}: $|X_1(u)| \le n_1-2$. \\
If $n_2=2$, then we also have $n_3=2$, and hence $G = K_{n_1} \cp K_{2} \cp K_{2} = K_{n_1} \cp C_4$. The assertion of the theorem then  follows by Proposition~\ref{prop:cp-cycle-by-complete graphs}. In the rest we may thus assume that $n_2\geq 3$.

Suppose on the contrary that $\mut(G)\geq n_1+n_2-1\geq n_1+2$. Then $|X_3(u)|\geq 3$. We distinguish two subcases.

\medskip\noindent
{\bf Case 2.1}: $n_3=2$. \\
Then the third coordinate of the vertices from $X_3(u)$ is $2$.
Let $z = (i,j,2)\in X_3(u)$. If no vertex from $X$ is adjacent to $z$, then $|X| \le |X_1(u)| + 2 \le n_1$ and we are done. Assume hence that $z'\in X$ is adjacent to $z$, so that $z'=(i',j,2)$ or $z'=(i,j',2)$. If $z'=(i',j,2)$, then the first coordinates of the vertices from $X$ are pairwise different, hence $|X|\le n_1$ in this subcase. If $z'=(i,j',2)$, then the vertices from $X_3(u)$ must have pairwise different second coordinates and also not equal to $1$, hence in this subcase $|X_3(u)| \le n_2 -1$ so that $|X| \le 1 + (n_1-2) + (n_2-1) = n_1 + n_2 - 2$.

\begin{figure}[ht!]
\begin{center}
\begin{tikzpicture}[scale=1.0,style=thick]
\tikzstyle{every node}=[draw=none,fill=none]
\def\vr{3pt} 

\begin{scope}[yshift = 0cm, xshift = 0cm]
\node [below=0.5mm ] at (0,0) {$u$};
\node [ ] at (0,2) {$\cdots$};
\node [ ] at (-2.5,2) {$\cdots$};
\node [ ] at (0,4) {$\cdots$};
\node [ ] at (2,6) {$\cdots$};
\node [ ] at (-1.5,6) {$\cdots$};
\node [ ] at (0.5,6) {$\cdots$};
\node [above=0.15mm ] at (-2.9,5.5){$w$};

\path (0,0) coordinate (x1);
\path (-3,2) coordinate (x2);
\path (-2,2) coordinate (x3);
\path (-1,2) coordinate (x4);
\path (1,2) coordinate (x5);
\path (2,2) coordinate (x6);
\path (3,2) coordinate (x7);
\path (-3,4) coordinate (x8);
\path (-2,4) coordinate (x9);
\path (-1,4) coordinate (x10);
\path (1,4) coordinate (x11);
\path (2,4) coordinate (x12);
\path (3,4) coordinate (x13);
\path (-3,6) coordinate (x14);
\path (-2,6) coordinate (x15);
\path (-1,6) coordinate (x16);
\path (0,6) coordinate (x17);
\path (1,6) coordinate (x18);
\path (3,6) coordinate (x19);

\draw (x1) --(x2);
\draw (x1) --(x3);
\draw (x1) --(x4);
\draw (x1) --(x5);
\draw (x1) --(x6);
\draw (x1) --(x7);
\draw (x14) --(x15);
\draw (-1.5,3.4) --(-2.5,2.6);
\draw (-1.3,3.4) --(0,2.6);
\draw (1.5,3.4) --(0.2,2.6);
\draw (1.7,3.4) --(3,2.6);
\draw (-1.5,5.1) --(-2.5,4.6);
\draw (-1.3,5.1) --(0,4.6);
\draw (1.5,5.1) --(0.2,4.6);
\draw (1.7,5.1) --(3,4.6);

\draw (x14) .. controls (-1.8,6.35).. (x16);
\draw (-2.5,2) ellipse (0.8 and 0.25);
\draw (0,2) ellipse (3.5 and 0.5);
\draw (0,4) ellipse (3.5 and 0.5);
\draw (0,6) ellipse (3.8 and 0.8);
\draw (-1,6) ellipse (2.4 and 0.7);
\draw (-1.5,6) ellipse (0.8 and 0.25);
\draw (0.5,6) ellipse (0.8 and 0.25);

\draw (x1)  [fill=black] circle (\vr);
\draw (x2)  [fill=black] circle (\vr);
\draw (x3)  [fill=black] circle (\vr);
\draw (x4)  [fill=white] circle (\vr);
\draw (x5)  [fill=white] circle (\vr);
\draw (x6)  [fill=white] circle (\vr);
\draw (x7)  [fill=white] circle (\vr);
\draw (x8)  [fill=white] circle (\vr);
\draw (x9)  [fill=white] circle (\vr);
\draw (x10) [fill=white] circle (\vr);
\draw (x11) [fill=white] circle (\vr);
\draw (x12) [fill=white] circle (\vr);
\draw (x13) [fill=white] circle (\vr);
\draw (x14) [fill=black] circle (\vr);
\draw (x15) [fill=black] circle (\vr);
\draw (x16) [fill=black] circle (\vr);
\draw (x17) [fill=black] circle (\vr);
\draw (x18) [fill=black] circle (\vr);
\draw (x19) [fill=white] circle (\vr);

 \draw[left] (x2)++(-0.5,0.0) node {$X_1(u)$};
 \draw[left] (x14)++(-0.8,0.0) node {$X_3(u)$};
 \draw[right] (x7)++(0.5,0.0) node {$N_1(u)$};
 \draw[right] (x13)++(0.5,0.0) node {$N_2(u)$};
 \draw[right] (x19)++(0.8,0.0) node {$N_3(u)$};
 \draw[above] (-1.5,6)++(0.0,0.15) node {{\small $W_1$}};
 \draw[above] (0.0,6)++(0.0,0.10) node {{\small $W_3$}};
\end{scope}
\end{tikzpicture}
\end{center}
\caption{Graph $G$ with a total mutual-visibility set $X=\{u\}\cup X_1(u)\cup X_3(u)$.}
\label{fig:W1W3}
\end{figure}
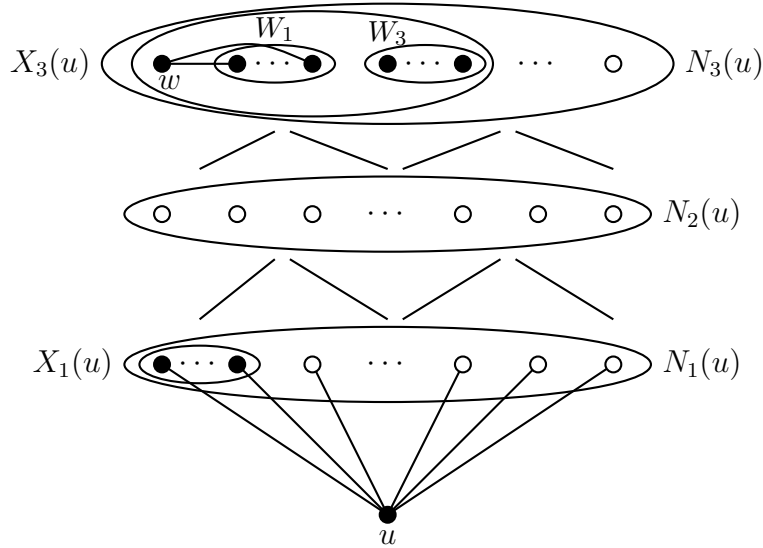

\medskip\noindent
{\bf Case 2.2}: $n_3=3$. \\
If the third coordinate of all the vertices of $X_3(u)$ is $2$, then by similar arguments to the case $n_3=2$ we get $\mut(G)\leq n_1+n_2-2$. Analogously, if the third coordinate of all the vertices of $X_3(u)$ is $3$, then we also conclude that $\mut(G)\leq n_1+n_3-2 \leq n_1+n_2-2$. Hence in the rest we may assume that in $X_3(u)$ there exist vertices $w=(i,j,2)$ and $w'=(i',j',3)$. Set $W_1 = X_3(u) \cap X_1(w)$ and $W_3 = X_3(u) \cap X_3(w)$, so that $X_3(u)=\{w\}\cup W_1\cup W_3$. Further, let $N_i(u)=\{v\in V(G):d_G(u,v)=i\}$, $i\in[3]$. See Fig.~\ref{fig:W1W3} for a schematic presentation of these sets.

Assume first that $W_1 = \emptyset$, that is, $X_3(u)=\{w\}\cup W_3$. Since $n_3=3$ and $w=(i,j,2)$,  the third coordinate of the vertices of $W_3$ is $3$. Furthermore, if the second coordinates of all the vertices from $W_3$ are equal to $j'$, the first coordinates of the vertices from $X$ are pairwise different. Then we have $|W_3|\leq  n_1-|X_1(u)|-2$ and hence
\begin{equation*}
\begin{aligned}
|X|&=|X_1(u)|+|X_3(u)|+1\\
&=|X_1(u)|+(|W_3|+1)+1\\
&\leq |X_1(u)|+(n_1-|X_1(u)|-2+1)+1=n_1.
\end{aligned}
\end{equation*}
If the first coordinates of all the vertices from $W_3$ are equal to $i'$, then since $j'\neq j\neq 1$, we have $|W_3|\leq n_2-2$. Hence in this case we get
\begin{equation*}
\begin{aligned}
|X|&=|X_1(u)|+|X_3(u)|+1\\
&=|X_1(u)|+(|W_3|+1)+1\\
&\leq n_1-2+(n_2-2+1)+1=n_1+n_2-2.
\end{aligned}
\end{equation*}

Assume second that $|W_1| > 0$. If the second coordinates of all the vertices of $W_1$ are $j$, then it follows that the third coordinates of the vertices from $W_1$ equal to $2$ and the third coordinates of the vertices from $W_3$ equal to $3$. Hence the first coordinates of the vertices from $X_1(u)\cup W_1$ are pairwise different, since $|W_3|\geq 1$, which then implies that $|W_1|\leq n_1-|X_1(u)|-2$.
 Assume further that the second coordinates of all the vertices from $W_3$ equal to $j'$. That means that the first coordinates of the vertices from $X$ are pairwise different, then we have $|W_3|\leq n_1-|X_1(u)|-|W_1|-2$.
 Hence
 \begin{equation*}
\begin{aligned}
|X|&=|X_1(u)|+|X_3(u)|+1\\
&=|X_1(u)|+(|W_1|+|W_3|+1)+1\\
&\leq |X_1(u)|+(|W_1|+n_1-|X_1(u)|-|W_1|-2+1)+1=n_1.
\end{aligned}
\end{equation*}
Assume that the first coordinates of all the vertices from $W_3$ are equal to $i'$. Since $j'\neq j\neq 1$, then we have $|W_3|\leq n_2-2$. Hence
\begin{equation*}
\begin{aligned}
|X|&=|X_1(u)|+|X_3(u)|+1\\
&=|X_1(u)|+(|W_1|+|W_3|+1)+1\\
&\leq |X_1(u)|+(n_1-|X_1(u)|-2+n_2-2+1)+1=n_1+n_2-2.
\end{aligned}
\end{equation*}
Similarly, we also get $|X|\leq n_1+n_2-2$ when the first coordinates of all the vertices of $W_1$ are $i$.

By similar arguments to the above, $\mut(G)\le n_1+n_2-2$ holds when $u$ and $u'$ differ in the second or in the third coordinate, in which case all the vertices from $X_1(u)$ differ in this coordinate. Hence in any case we have $\mut(G)\le n_1+n_2-2$ and we can conclude that $\mut(G)=n_1+n_2-2$.
\qed

\begin{theorem}
	\label{thm:three-complete graphs-large-cases}
	If $n_1\ge n_2\ge n_3\ge 4$, then $\mut(K_{n_1}\cp K_{n_2}\cp K_{n_3}) = n_1+n_2+n_3-6$.
\end{theorem}

\proof
To prove this result we will apply Proposition~\ref{prop:reduction}. More precisely, setting $H = K_{n_1,n_2,n_3}$ we are going to prove that the largest set of triangles in $H$ such that no two triangles from the set intersect in a single vertex has cardinality $n_1 + n_2 + n_3 - 6$.

We use the same notation as in the proof of Proposition~\ref{prop:reduction}. That is, let $I_1, I_2, I_3$ be the partition classes of $H$ with $I_i = \{u_{i,j}:\ j \in [n_i]\}$ for each $i \in [3]$. A $3$-clique (triangle) induced by the vertices $u,v,z$ is  denoted by $uvz$. Consider first the following sets of triangles in $H$
$$\cK_1=\{ u_{1,j}u_{2,1}u_{3,1} :  3\leq j \leq n_1 \},  \qquad
\cK_2=\{ u_{1,1}u_{2,j}u_{3,2} :  3\leq j \leq n_2 \}, $$
$$\cK_3=\{ u_{1,2}u_{2,2}u_{3,j} :  3\leq j \leq n_3 \},  $$
and set $\cK= \cK_1 \cup \cK_2 \cup \cK_3$. It is clear that $|\cK|= n_1+n_2+n_3-6$. Moreover, for any two triangles $t, t' \in \cK$, either $t$ and $t'$ are from the same set $\cK_i$ and $|t\cap t'|=2$, or they are from different triangle sets and do not share a vertex. In either case, $|t \cap t'| \neq 1$ as required. By Proposition~\ref{prop:reduction} we thus have $\mut(K_{n_1}\cp K_{n_2}\cp K_{n_3}) \ge n_1+n_2+n_3-6$.
\medskip

To prove the reverse inequality, let $\cK$ be a triangle set in $H$ of maximum cardinality such that $|t\cap t'|\neq 1$ holds for every pair of triangles $t,t'$ from $\cK$. Given $\cK$, we say that  $uv$ is a \emph{base edge} of the triangle $uvz \in \cK$, if $uv$ is incident to at least two triangles from $\cK$.
\paragraph{Claim A.} Every triangle $t \in \cK$ has at most one base edge.\\
\emph{Proof. } Suppose that $uv$ and $vz$ are two base edges of $uvz \in \cK$. Then, there exists a vertex $z_1 \neq z$ with $uvz_1 \in \cK$, and also a vertex $u_1\neq u$ with $u_1vz \in \cK$. The triangles $t_1=uvz_1$ and $t_2= u_1vz$ share the vertex $v$. As $|t_1 \cap t_2|=1$ is not possible, there must be another common vertex in $t_1 \cap t_2$. As $u\neq u_1$, $z\neq z_1$, and $u \neq z$ are supposed, the only remaining possibility is $u_1=z_1$. However, under the assumption $u_1=z_1$, vertices  $u$, $v$, $z$, and $u_1$ form a  $4$-clique  that is impossible in the $3$-partite graph $H$.
\smallqed
\medskip

By Claim~A, the triangles in $\cK$ can be partitioned into classes $\cK_1, \dots , \cK_s$ such that, for each $i \in [s]$, the set $\cK_i$ either contains all triangles from $\cK$ which are incident to a fixed base edge, or $\cK_i$ contains just one triangle without a base edge. Note that, by Claim~A, the partition is unique. Let $V_i$ denote the set of vertices covered by the triangles in $\cK_i$, for every $i \in [s]$.
\paragraph{Claim B.} The sets $V_1, \dots ,V_s$ are pairwise vertex-disjoint, and $|\cK|= |\bigcup_{i=1}^{s} V_i| -2s \leq n_1+n_2+n_3-2s$ holds.  \\
\emph{Proof. } We first prove that $V_1, \dots ,V_s$ are pairwise vertex-disjoint. If $\cK_i$ contains a triangle $t$ without a base edge, then $|t\cap t'| <2$ holds for every $t' \in \cK \setminus \{t\}$. By the assumption $|t \cap t'| \neq 1$, we may infer that $t$ is vertex-disjoint from every other triangle in $\cK$.
Now, consider a triangle $t=uvz$ with a base edge $uv$ and the class $\cK_i$ that contains $t$. Suppose that $z$ is also  incident to another triangle $t'=zxy$. As $|t \cap t'| \neq 1$, the two triangles share a vertex different from $z$. It follows then that $t$ contains a base edge different from $uv$ that contradicts Claim~A. We may conclude that $z$ belongs to only one class $V_i$. A similar argument shows that the same is true for the vertices $u$ and $v$.

Since the triangles in $\cK_i$ share a base edge or $\cK_i$ contains only one triangle, $|\cK_i|= |V_i|-2$ holds for every $i \in [s]$. As $V_1, \dots , V_s$ are pairwise vertex-disjoint, we may conclude  $|\cK|= |\bigcup_{i=1}^{s} V_i| -2s \leq n_1+n_2+n_3-2s$ as stated.
\smallqed
\paragraph{Claim C.} If $s=1$, then $|\cK|\leq  n_1$, and if $s=2$, then  $|\cK| \leq  n_1+n_2-2$.\\
\emph{Proof. } If $s=1$, then all triangles in $\cK$ has the same base edge $uv$ and a further vertex which belongs to the third partite class of $H$. Hence, the number of triangles in $\cK$ is at most $n_1$ and $|\cK| \leq n_1+n_2-2$ is true.

If $s=2$, we have two sets $\cK_1$ and $\cK_2$ with base edges $uv$ and $u'v'$. If $u,u' \in I_i$ and $v,v' \in I_j$, then the remaining vertices in $V_1 \cup V_2 \setminus \{u,u',v,v' \}$ all belong to the third partite class $I_\ell$. By Claim~B, $V_1$ and $V_2$ are disjoint sets and therefore, $|V_1 \cup V_2| \leq |I_\ell| +4$ and $|\cK| \leq |I_\ell|=n_\ell$. In the other case, the base edges $uv$ and $u'v'$ contain at least one vertex from each partition class of $H$. Suppose that $u, u' \in I_i$, $v \in I_j$,  $v' \in I_\ell$, and that $uv$ and $u'v'$ are the base edges in $\cK_1$ and $\cK_2$, respectively. It follows from Claim~B that $V_1 \setminus \{u,v\} \subseteq I_\ell \setminus \{v'\}$ and $|\cK_1| \leq n_\ell-1$. Analogously, $|\cK_2| \leq n_j-1$. We may conclude again that
$$|\cK| \leq n_\ell +n_j -2 \leq n_1+n_2 -2$$
which finishes the proof of the claim.
\smallqed
\medskip

By Claim~B, we have $|\cK| \leq n_1+n_2+n_3-6$ whenever $s \ge 3$.  By Claim~C, we get $|\cK| \leq n_1+n_2 -2 $ if $s=2$. As $n_3 \ge 4$ is supposed, it also implies
$|\cK| \leq n_1+n_2+n_3-6$. Similarly, the case of $s=1$ gives $|\cK| \leq n_1< n_1+n_2+n_3-6$. We may then infer $|\cK| \leq n_1+n_2+n_3-6$. In view of Proposition~\ref{prop:reduction}, we conclude that $\mut(K_{n_1}\cp K_{n_2}\cp K_{n_3}) \le n_1+n_2+n_3-6$.
\qed

Theorem~\ref{thm:three-complete graphs} follows by combining Theorems~\ref{thm:three-complete graphs-small-cases} and~\ref{thm:three-complete graphs-large-cases}.

\begin{remark}
Using the method from the proof of Theorem~\ref{thm:three-complete graphs-large-cases}, it is possible to give a shorter proof of Theorem~\ref{thm:three-complete graphs-small-cases}. However, we decided to include the present proof because it could provide a different technique to handle higher dimensional Hamming graphs. 
\end{remark}

\section{Proof of Theorem~\ref{thm:upper}}
\label{sec:upper}

Recall that for Theorem~\ref{thm:upper}$(i)$ we assume $r \ge 3$ and $N = n_1+\cdots + n_r$. The main goal is to prove that $\mut(K_{n_1} \cp \cdots \cp K_{n_r}) = \cO(N^{r-2})$, but along the way we will also determine the constant from the statement of the theorem. We are going to show that, for each $r \ge 3$ and each $r$-dimensional Hamming graph $K_{n_1} \cp \cdots \cp K_{n_r}$, the constant $c_r = \frac{6}{r!}$ satisfies
\begin{equation} \label{eq:upper}
\mut(K_{n_1} \cp \cdots \cp K_{n_r}) \leq  c_r N^{r-2}\, .
\end{equation}
By Theorem~\ref{thm:three-complete graphs}, the statement~(\ref{eq:upper}) is true for $r=3$ with the constant $c_3=1$. We then proceed by induction on $r$.

We may suppose that $n_1 \ge \cdots \ge n_r$ and therefore, $n_r \leq \frac{N}{r}$ holds. Let $H$ denote the product $K_{n_1} \cp \cdots \cp K_{n_r}$ and let $X$ be a maximum total mutual-visibility set in $H$. For an integer $j \in [n_r]$, we define the following set of vertices in $H$:
$$ L(j)= \{(x_1, \dots , x_{r-1}, j):  x_i \in [n_i] \mbox{ for  } i \in [r-1]  \} .  $$
The subgraph of $H$ induced by $L(j)$ is isomorphic to the $(r-1)$-dimensional Hamming graph $H'= K_{n_1}\cp \cdots \cp K_{n_{r-1}}$. As it is a convex subgraph of $H$, the set $X \cap L(j)$ contains no two vertices at distance $2$ apart. Hence, the removal of the fixed last coordinate $j$ transforms $L(j) \cap X$ into a total mutual-visibility set $X'$ in $H'$. Clearly, $|X'| \leq \mut(H')$ and $|X'|=|X \cap L(j)| $ hold. By the induction hypothesis,
\begin{equation} \label{eq:A}
|X'| \leq \mut(H') \leq c_{r-1} (N-n_r)^{r-3} <  c_{r-1}\, N^{r-3}.
\end{equation}
Equivalently, every set $ L(j)$ contains less than $c_{r-1}\, N^{r-3}$ vertices from $X$. On the other hand, every vertex $v \in X$ belongs to exactly one set $L(j)$ and it follows that
\begin{equation} \label{eq:B}
	|X| < n_r\, c_{r-1}\, N^{r-3} \leq \frac{N}{r} \, c_{r-1}\, N^{r-3} = \frac{c_{r-1}}{r}\, N^{r-2} .
\end{equation}
This proves~(\ref{eq:upper}) for the $r$-dimensional Hamming graphs with the constant $c_r= c_{r-1}/r$. Starting with the constant $c_3=1$, we infer that  the upper bound~(\ref{eq:upper}) holds for every $r \ge 3$ with the constant $c_r=\frac{6}{r! }$.
\medskip

The assertion $(ii)$ of Theorem~\ref{thm:upper} is true for $r=3$ by Theorem~\ref{thm:three-complete graphs}. We then proceed by induction on $r$. The formula can be proved along the same lines as the inequality in $(i)$. We set $H=K_s^{\cp, r}$,  $n_1=\dots =n_r=s$, and $N=rs$.
 Rewriting the inequalities $(\ref{eq:A})$ and $(\ref{eq:B})$ according to the hypothesis and using $N-n_r=(r-1)s$, we get
\begin{equation*} 
	|X'| \leq \mut(H') \leq \left(3\prod_{i=3}^{r-1} (i-1)^{i-3} \right) ((r-1)s)^{r-3} = \left( 3\prod_{i=3}^{r} (i-1)^{i-3} \right)\, s^{r-3},
\end{equation*}
and we can conclude
\begin{equation*} 
	|X| =\mut(H) \leq  \left( 3\prod_{i=3}^{r} (i-1)^{i-3} \right) \, s^{r-2}
\end{equation*}
which proves Theorem~\ref{thm:upper}.

\section{Proof of Theorem~\ref{thm:lower}}
\label{sec:balanced}

Theorem~\ref{thm:lower} asserts that for every integer $r\ge 3$ it holds that $\mut(K_s^{\cp, r}) = \Theta(s^{r-2})$. For $r \ge 3$, Theorem~\ref{thm:upper}$(ii)$ directly implies $\mut(K_s^{\cp, r}) = \cO(s^{r-2}).$ For the lower bound, we give a probabilistic proof based on a similar idea as the proof in~\cite[Section~4]{BES} for a famous hypergraph Tur\' an-problem of Brown, Erd\H{o}s, and S\' os.

Let $r \ge 3$ and  $H=K_s^{\cp, r}$. From $H$, we choose each vertex with probability $p=\frac{2}{r(r-1)s^2}$, independently of the decisions made for other vertices. This way, we obtain a set $S \subseteq V(H)$. The expected value of the size of $S$ is
$$E(|S|) = s^r p= \frac{2}{r(r-1)} s ^{r-2}\, .$$
We say that a set $\{u,v\}$ of two vertices from $S$ is a \emph{bad pair} in $S$, if $d_H(u,v)=2$. Let $B$ denote the set of all bad pairs that are present in $S$. To get the pairs of vertices at distance $2$ apart in $H$, we may first choose the two entries where they differ, fix the coordinates for these two entries appropriately, and fix the remaining coordinates arbitrarily. As the vertices of $S$ were selected randomly with probability $p$, we may estimate the size of $B$ in the following way:
\begin{equation*}
\begin{aligned}
 E(|B|) & = {r \choose 2} \frac{s^2 (s-1)^2 }{2} s^{r-2}\, p^2\\
       & \leq \frac{r(r-1) s^4}{4} \, \frac{2}{r(r-1)s^2}\,  s^{r-2} \, p\\
       & = \frac{1}{2}\, s^r \, p =  \frac{1}{2}\, E(|S|).
\end{aligned}
\end{equation*}

For a set $S \subseteq V(H)$, we remove one vertex from each bad pair. The obtained set $S^*$ contains no bad pairs, that is, $S^* $ is a total mutual-visibility set in $H$. Moreover, we have
$$E(|S^*|) \geq E(|S|)- E(|B|) \geq \frac{1}{2}\, E(|S|)= \frac{s^{r-2}}{ r(r-1) }.$$
There exists at least one set $S$ that results in a total mutual visibility set $X$ with  $|X| \geq E(|S^*|) $ after removing one vertex from each bad pair. We may therefore infer
$$\mut(K_s^{\cp, r}) \geq \frac{1}{ r(r-1) }\, s^{r-2} .$$
Together with the upper bound, this implies $\mut(K_s^{\cp, r}) = \Theta(s^{r-2})$.

\begin{remark}
If we consider $r$-dimensional Hamming graphs in general, the statement analogous to that in Theorem~\ref{thm:lower} is not valid. Indeed, assume that, for a fixed $r \ge 3$, there exists an absolute constant $c=c(r)$ such that
$$  \mut(K_{n_1} \cp \cdots \cp K_{n_r} ) \geq c\, N^{r-2}$$
holds for every Hamming graph with $2\leq n_r \leq \cdots \leq n_1$ and $N=\sum_{i=1}^{r}n_i$. By setting $n_2= \cdots =n_{r}=2$ and choosing an integer $n_1 > 2 \sqrt[r-3]{4/c}$, the obtained $r$-dimensional Hamming graph $H$ gives the contradiction
$$\mut(H) \leq |V(H)| = 2^{r-1}\, n_1 = \frac{2^{r-1}}{n_1^{r-3}}\, n_1^{r-2} < c\, n_1^{r-2} < c\, N^{r-2} \leq \mut(H).$$
\end{remark}

\section{A Tur\'an-type problem}
\label{sec:conclude}

In this section, we show that our main results can be reformulated in the language of Tur\'an-type problems on hypergraphs.

A \emph{hypergraph} $H=(V,E)$ is a set system over the vertex set $V$. More precisely, every (hyper)edge $e \in E$ is a nonempty subset of $V$. A hypergraph $H'=(V', E')$ is a \emph{subhypergraph} of $H=(V,E)$ if both $V' \subseteq V$ and $E' \subseteq E$ hold. A hypergraph $H$ is \emph{$r$-uniform} if every $e \in E$ contains exactly $r$ vertices. Note that $r$-uniform hypergraphs are often called {\em $r$-graphs}. Then, $2$-uniform hypergraphs correspond to simple graphs. The \emph{complete $r$-partite $r$-graph $\cK_{n_1, \dots, n_r}^{(r)}$} is the $r$-uniform hypergraph on the vertex set $V=V_1 \cup \cdots \cup V_r$, where the partite classes $V_1, \dots, V_r$ are pairwise disjoint and $|V_i|=n_i$ holds for every $i \in [r]$, and moreover, the edge set is defined as
$E=\{e \subseteq V: |e \cap V_i|= 1 \mbox{ for every } i \in [r] \}$. Thus, $\cK_{n_1, \dots, n_r}^{(r)}$ contains $\Pi_{i=1}^r n_i$ edges.

The basic example for a hypergraph Tur\'an problem takes an $n$-element vertex set $V$ and asks for the maximum number of edges in a  $r$-uniform hypergraph $H=(V,E)$ that contains no subhypergraph isomorphic to a given ($r$-uniform) hypergraph $F$. This maximum number is denoted by $\ex_r(n, F)$. Remark that most of the Tur\'an-type hypergraph problems considered are notoriously hard. Even the tight asymptotics or the exact order of magnitude for $\ex_r(n,F)$ may be hard to identify as $n \to \infty.$ For more details on the subject we refer the reader to the book~\cite{gerbner-2019} and the survey~\cite{keevash-2011}.

We may also consider the version of the problem, where the edges must be selected from the complete $r$-uniform $r$-graph $\cK_{n_1, \dots, n_r}^{(r)}$ such that the obtained hypergraph does not contain a subhypergraph isomorphic to a given $F$. Under this condition, the maximum number of edges will be denoted by $\ex(\cK_{n_1, \dots, n_r}^{(r)}, F)$.

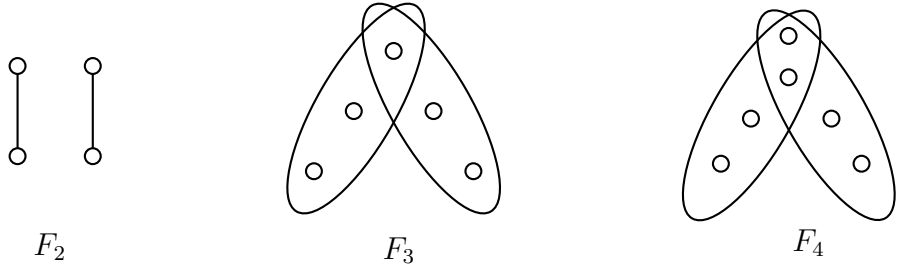
\begin{figure}[ht!]
\begin{center}
\begin{tikzpicture}[scale=1.0,style=thick]
\tikzstyle{every node}=[draw=none,fill=none]
\def\vr{3pt} 

\begin{scope}[yshift = 0cm, xshift = 0cm]
\path (0,-0.1) coordinate (x1);
\path (0.53,-0.9) coordinate (x2);
\path (1.06,-1.7) coordinate (x3);
\path (-0.53,-0.9) coordinate (x4);
\path (-1.06,-1.7) coordinate (x5);
\path (-4,-0.3) coordinate (x6);
\path (-4,-1.5) coordinate (x7);
\path (-5,-0.3) coordinate (x8);
\path (-5,-1.5) coordinate (x9);
\path (5.25,0.1) coordinate (x10);
\path (5.25,-0.45) coordinate (x11);
\path (5.82,-1) coordinate (x12);
\path (6.22,-1.6) coordinate (x13);
\path (4.75,-1.0) coordinate (x14);
\path (4.35,-1.6) coordinate (x15);

\draw (x6) --(x7);
\draw (x8) --(x9);

\draw[rotate=30] (0,-1.00) ellipse (15pt and 45pt);
\draw[rotate=-30] (0,-1.00) ellipse (15pt and 45pt);
\draw[rotate=30] (4.5,-3.7) ellipse (15pt and 45pt);
\draw[rotate=-30] (4.6,1.55) ellipse (15pt and 45pt);

\draw (x1)  [fill=white] circle (\vr);
\draw (x2)  [fill=white] circle (\vr);
\draw (x3)  [fill=white] circle (\vr);
\draw (x4)  [fill=white] circle (\vr);
\draw (x5)  [fill=white] circle (\vr);
\draw (x6)  [fill=white] circle (\vr);
\draw (x7)  [fill=white] circle (\vr);
\draw (x8)  [fill=white] circle (\vr);
\draw (x9)  [fill=white] circle (\vr);
\draw (x10) [fill=white] circle (\vr);
\draw (x11) [fill=white] circle (\vr);
\draw (x12) [fill=white] circle (\vr);
\draw (x13) [fill=white] circle (\vr);
\draw (x14) [fill=white] circle (\vr);
\draw (x15) [fill=white] circle (\vr);

 \draw[left] (x9)++(0.8,-1.2) node {$F_2$};
 \draw[left] (x5)++(1.5,-1.05) node {$F_3$};
 \draw[right] (x15)++(0.8,-1.05) node {$F_4$};
 
\end{scope}
\end{tikzpicture}
\end{center}
\caption{The $2$-uniform, $3$-uniform and $4$-uniform forbidden subhypergraphs in our problem.}
\label{fig:forbidden subhypergraphs}
\end{figure}

For $r \ge 2$, let $F_r$ denote the $r$-uniform hypergraph on $r+2$ vertices that contains two edges $f_1$ and $f_2$ with $|f_1 \cap f_2|= r-2$, see Fig.~\ref{fig:forbidden subhypergraphs}. 

By Proposition~\ref{prop:reduction}, the maximum size of a total mutual visibility set in the Hamming graph $K_{n_1} \cp \cdots \cp K_{n_r}$ equals the maximum number of $r$-cliques in the graph $K_{n_1, \dots, n_r}$ such that no two of them intersect in $r-2$ vertices. The latter problem can be expressed by taking the vertex sets of the maximum cliques in $K_{n_1, \dots, n_r}$ as edges in an $r$-uniform $r$-graph and forbidding the subhypergraph $F_r$. Therefore, we may conclude 
$$\mut(K_{n_1} \cp \cdots \cp K_{n_r})= \ex(\cK_{n_1, \dots, n_r}^{(r)}, F_r)$$
and then, Theorems~\ref{thm:three-complete graphs}, \ref{thm:upper}, and \ref{thm:lower}  can be reformulated as follows.

\begin{proposition}
	If $n_1\geq n_2\geq n_3\ge 2$ and $n = n_1 + n_2 + n_3$, then
	\[
	\ex(\cK_{n_1,n_2,n_3}^{(3)}, F_3) = \left.
	\begin{cases}
		n - 4; & n_3 = 2, \\
		n - 5; & n_3 = 3, \\
		n - 6; & n_3\ge 4\,.
	\end{cases}
	\right.
	\]
\end{proposition}
\begin{proposition}
\begin{itemize}
	\item[$(i)$] If $r \ge 3$ and $n$ denotes $\sum_{i=1}^r n_i$, then $$\ex(\cK_{n_1,\dots , n_r}^{(r)}, F_r) = \cO(n^{r-2}).$$
	\item[$(ii)$] For every integer $r \ge 3$, it holds that  $$\ex(\cK_{s,\dots , s}^{(r)}, F_r) = \Theta(s^{r-2}).$$
\end{itemize}
\end{proposition}

The famous problem of Brown, Erd\H{o}s, and S\'os from~\cite{BES} asks for the maximum number of edges in an $r$-uniform hypergraph of order $n$ when all subhypergraphs with $v$ vertices and $e$ edges are forbidden. This maximum is denoted by $f^{(r)}(n, v, e)-1$. Our problem differs from this famous one in two main aspects and  consequently, neither lower nor upper bounds on $f^{(r)}(n, r+2, 2)-1$ can be applied directly to $\ex(\cK_{n_1,\dots , n_r}^{(r)}, F_r)$. First, when $f^{(r)}(n, r+2, 2)-1$ is counted, $r$-edges intersecting in $r-1$ vertices are also forbidden unlike to our problem setting. Second, in our problem, the edges of the extremal hypergraph must be selected from $\cK_{n_1, \dots, n_r}$, while the problem discussed in~\cite{BES} has no such a restriction.

\section*{Acknowledgements}

This work was supported by the Slovenian Research Agency (ARIS) under the grants P1-0297, J1-2452, and N1-0285.


\end{document}